\documentclass{amsart}

\usepackage{amssymb}
\usepackage{amsmath}
\usepackage{latexsym,color}
\usepackage{amsthm,amsfonts, tikz, caption}
\usepackage{srcltx}

\usepackage{titlesec}

\titleformat{\subsection}
{\normalfont\fontsize{13}{13}\itshape}{\thesubsection}{1em}{}

\newcommand{\cE}{\ensuremath{\mathcal E}}

\newcommand{\var}{{\rm Var}} 

\newcommand{\bbR}{{\ensuremath{\mathbb R}} }

\newcommand{\bbZ}{{\ensuremath{\mathbb Z}} } 

\newcommand{\Z}{\mathbb Z}

\newcommand{\p}{{+}}
\newcommand{\m}{{-}}

\newcounter{mycount}

\title{The Life and Mathematical Legacy of Thomas M. Liggett}

\author[D. Aldous]{
  David Aldous}
  \address{
    David Aldous is a Professor Emeritus of Statistics at the
    University of California, Berkeley. His email address is aldousdj@berkeley.edu.
    }
  
 \author[P. Caputo]{ Pietro Caputo}
  \address{
  	Pietro Caputo is a Professor of Mathematics at Roma Tre University. His email address is  pietro.caputo@uniroma3.it.
   }

  \author[R. Durrett]{Rick Durrett}
\address{
	Rick Durrett is a Professor of Mathematics at
	Duke University. His email address is rtd@math.duke.edu.
}

  \author[P. Jung]{Paul Jung}
\address{
Paul Jung is an Associate Professor of Mathematical Sciences at the
Korea Advanced Institute of Science and Technology. His email address is pauljung@kaist.ac.kr.
}

  \author[A. E. Holroyd]{Alexander E. Holroyd} 
\address{  Alexander E. Holroyd  is a Professor of Mathematics at
	University of Bristol. His email address is a.e.holroyd@bristol.ac.uk.
}

 \author[A. L. Puha]{ Amber L. Puha}
\address{Amber L. Puha is a Professor of Mathematics at California State University, San Marcos. Her email address is apuha@csusm.edu.
}

\begin{document}

\begin{abstract}
	Thomas Milton Liggett was a world renowned UCLA probabilist, famous for his monograph {\it Interacting Particle Systems}. He passed away peacefully on May 12, 2020. This is a  perspective article in memory of both Tom Liggett the person and Tom Liggett the mathematician.
\end{abstract}

\maketitle

\begin{figure}[h]
	\centering
	\includegraphics[width=2.7in]{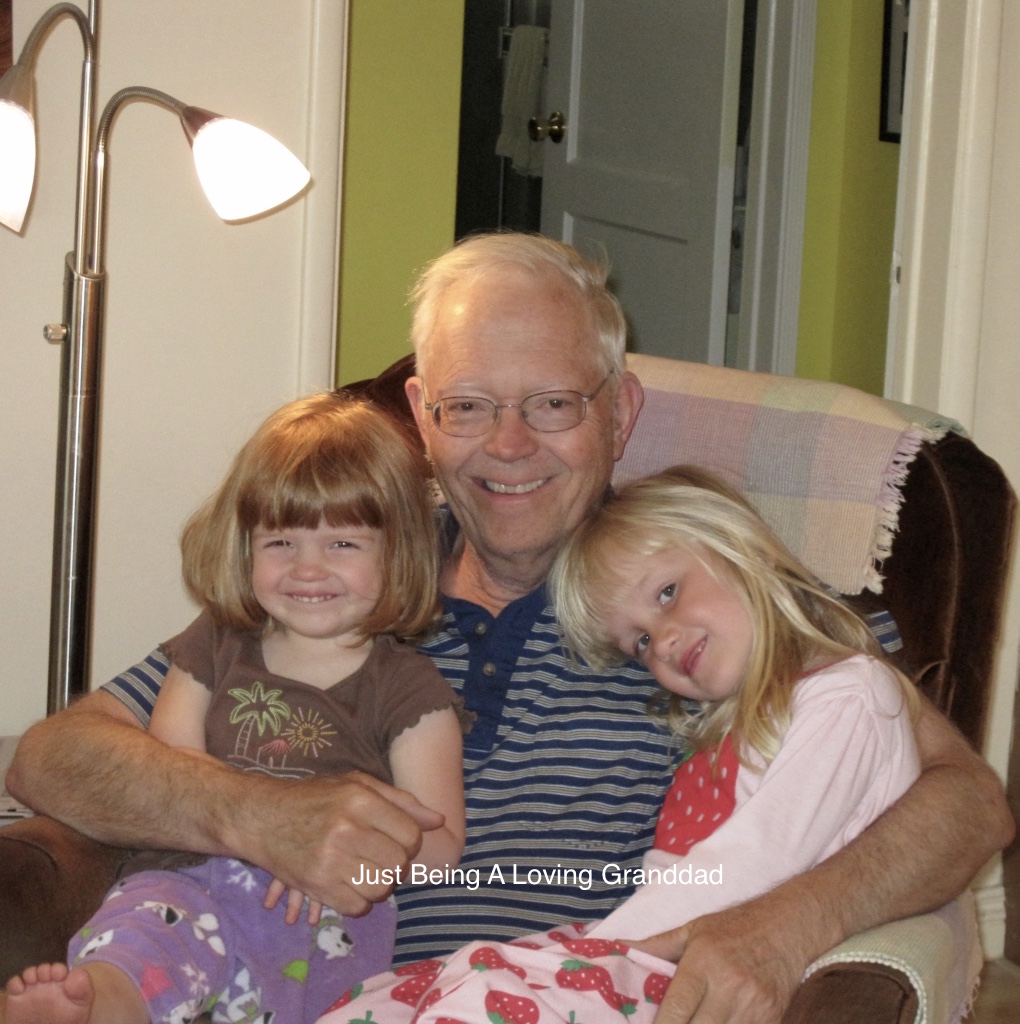}
	\caption*{Tom Liggett with his two granddaughters}
	\label{fig:Reinforcement}
\end{figure}

\begin{figure}[h]
	\centering
	\includegraphics[angle=90, width=2.7in]{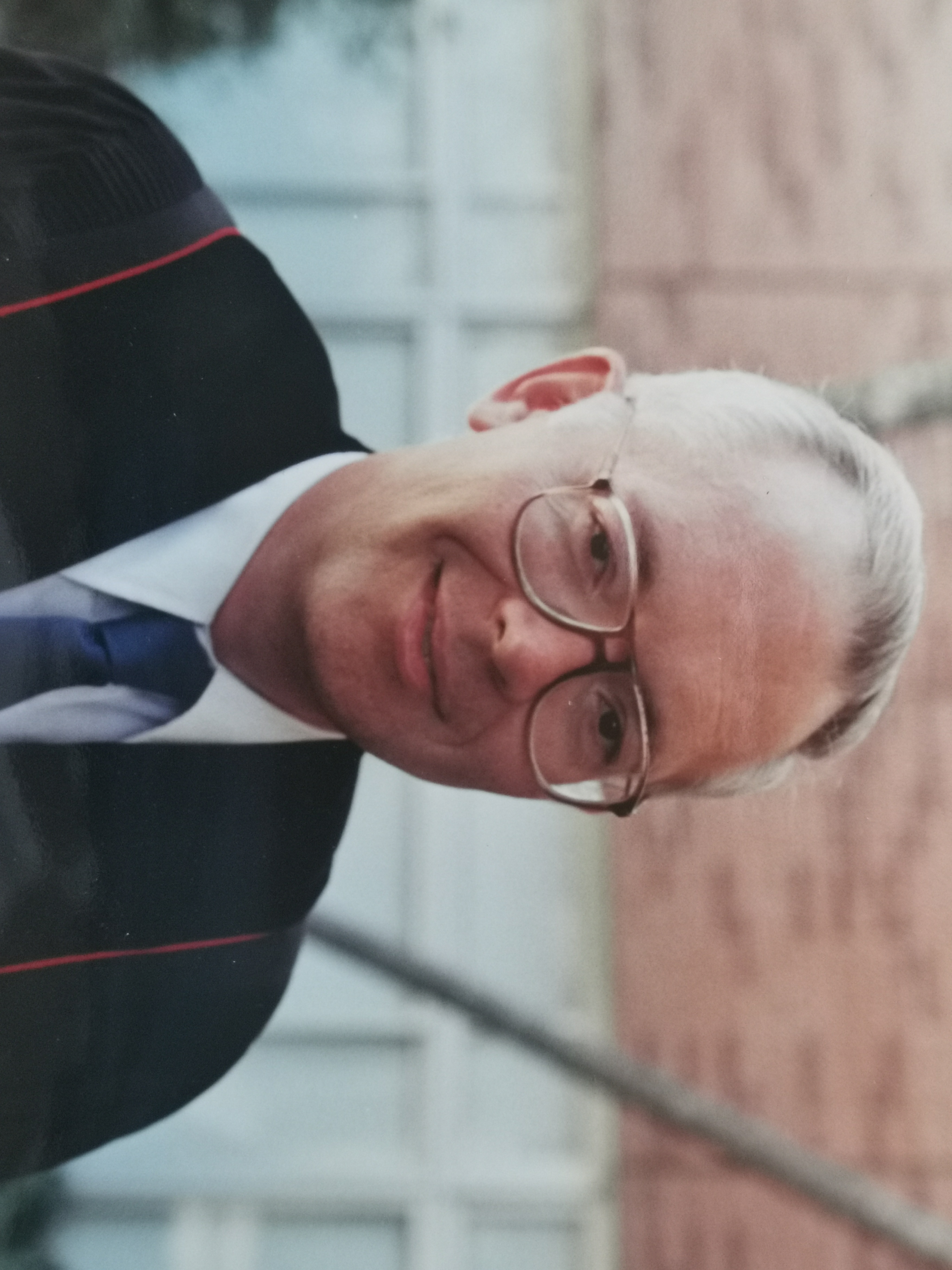}
	\caption*{Tom Liggett in 1989 at the graduation ceremony of his Ph.D. student Dayue Chen}
	\label{fig:Reinforcement}
\end{figure}

\section*{}

	\begin{center}
		{\large\bf Remembering Tom}
		
		\vspace{2mm}
		
		Paul Jung
	\end{center}
	
	Thomas Milton Liggett, a world renowned UCLA probabilist, passed away peacefully on May 12, 2020. He is survived by his wife, Christina Liggett, his son and daughter, Tim Liggett and Amy Liggett, and two granddaughters, Amanda Liggett and Jenna Liggett.
	
	Tom, as he was known affectionately to friends and colleagues, was born March 29, 1944 into a missionary family and spent much of his childhood in Argentina. He attended Oberlin College, and after graduating in 1965, he began his Ph.D. studies in  mathematics at Stanford University. His eulogy from the math department webpage at UCLA, written by colleagues Marek Biskup and Roberto Schonmann, reads ``Liggett's interest in probability was spawned by interactions with Samuel Goldberg during  his  undergraduate  time  at  Oberlin  College.   He continued  to  a  PhD  program at Stanford where he was further influenced by lectures of probabilist Kai Lai Chung. He ultimately wrote his PhD under the supervision of Samuel Karlin whom he found a better match personally,  but  he  did  not  find  probability  research  exciting  and  was even readying himself for a career of a liberal-arts college lecturer, rather than a research mathematician. His advisor urged him to at least call UCLA Math and ask for job application forms.  To his surprise, the letter he received in return contained a job offer and this is how he ended up moving to Southern California in 1969.''
	
	Upon arriving at UCLA in 1969, Tom met his future wife Christina Marie Goodale, who was working as an administrator in the math department. Their courtship began in 1971 and continued through Tom's first sabbatical visiting Jacques Neveu at Paris VI. More than one year later, along with forty-four letters (hand-written snail mail of course) and  an academic approval from UCLA (to avoid nepotism), they were married in August of 1972. Together Tom and Chris had two children.
	Their first child, Tim, enjoyed math and science just like Tom and became a high school physics teacher  while their daughter, Amy, emulated her paternal grandparents to become a minister in the church. Those of us who knew Tom well, knew him to be both a serious mathematician and a devoted father. A humorous account of growing up in the Liggett family, was given in a speech by Tim Liggett at Tom's 75th birthday conference which one can listen to by visiting celebratio.org/Liggett\_TM.
	
	Regarding his mathematics, Tom's name was practically synonymous with his famous monograph {\it Interacting Particle Systems} \cite{MR776231}. Many of his most important results were related to this subfield of probability theory, starting with his most cited research paper which is not a result in probability, but rather a result in functional analysis which extends the Hille-Yosida theorem to nonlinear generators \cite{MR287357}. This was a stepping stone which led Tom to prove a general existence theorem for particle systems on infinite graphs \cite{MR309218}. This existence result was a fundamental step in bringing the then burgeoning field to its modern form. A fuller account of Tom's enormous impact in the field of interacting particle systems is given later in this article by Rick Durrett. Let us just highlight one of Tom's big results in this field from \cite{MR806224} which improves upon Kingman's well-known subadditive ergodic theorem. These days, this versatile tool is often called the Kingman-Liggett subadditive ergodic theorem (see for instance, Fristedt and Gray's textbook {\it A Modern Approach to Probability Theory}).
	
	While particle systems were Tom's bread and butter, he also proved key results in several other related areas of probability theory. In \cite{MR2629990}, Tom was part of a team which proved the longstanding Aldous Spectral Gap Conjecture in the field of Mixing Times. An account of how this came to fruition is given later in this article by Pietro Caputo, with an addendum by David Aldous. A recurring theme of Tom's work was his keen problem solving ability. One recounting of this ``superpower'' of his, is given later in this article by Ander Holroyd in the context of their joint work \cite{MR3570073} on finitely dependent processes.   One of Tom's several papers on negative dependence \cite{MR2476782}, while published over a decade ago, is still considered one of the most important and relevant papers in this area. An account of Tom's influence in the area of negative dependence, by Robin Pemantle, is due to appear in a volume of Celebratio Mathematica in honor of Tom Liggett (celebratio.org/Liggett\_TM). Also in that same volume, one can find Tom's final (unpublished) manuscript concerning a family of random variables related to negative dependence, as well as a number of tributes to him written by fellow mathematicians and friends.
	
	In addition to his research Tom cared deeply about teaching mathematics. Like his books and research papers, his teaching was crystal clear. When explaining a proof, he effortlessly pointed out how each assumption in the theorem played its role. Unusual for a top-level mathematician like Tom, he took great care not only for graduate level teaching, but also for the teaching of introductory undergraduate level courses. For example, Tom once requested to teach a three-quarter consecutive sequence of introductory calculus just so that the students could have a continuous learning experience at this foundational level.
	
	It was Tom's motivational and pedagogical lecture style that convinced me to pivot from studying analysis to probability, under his guidance, for my Ph.D. dissertation. It was a decision I have never once regretted. Tom's guidance was ever thoughtful-- knowing exactly when to nudge, when to gently admonish, and when and how to encourage. His advising was selfless, and he remained an ever willing and open mentor for almost two decades after his advising duties were over. 
	
	Tom's mathematical brilliance was recognized by many awards which include, among others, a membership in the NAS and AAAS, an invitation to speak at the 1986 ICM, Sloan and Guggenheim Fellowships, and being named a Fellow of both the AMS and IMS. Tom's career was celebrated, on the occasion of his 65th birthday in 2009, at a conference in Beijing in his honor. For his 75th birthday, there was a conference in March 2019 in his honor at IPAM on the UCLA campus. Unfortunately, Tom developed a severe case of pneumonia right before the conference and was unable to attend. The toll from the pneumonia eventually led to hospice care until he left us in May of 2020. While his person may be gone, his spirit, legacy and outsized influence live on.

\section*{}

	\begin{center}
		{\large\bf Thomas M. Liggett the Thesis Advisor}
		
		\vspace{2mm}
		
		Amber L. Puha
	\end{center}

	Thomas Milton Liggett was my thesis advsior.  Over his 42-years as a faculty member in the UCLA Department of Mathematics, Tom had 9 such students.  {Norman Matloff (1975, UC Davis), Diane Schwartz (1975, CSU Northridge), Enrique Andjel (1981, U. Provence, Marseille, FR),  Dayue Chen (1989, Peking University), Xijian Liu (1991, US Census Bureau), Shirin Handjani, 1993 (CCR, San Diego), myself (1998, CSU San Marcos), Paul Jung (2003, KAIST) and Alexander Vandenberg-Rodes (2011).} I was lucky number 7.  To some this might seem a small number.  He once told me that he found such interactions unnatural.  Maybe so.  But, he was a perfect thesis advisor for me, and for 8 others.  He pushed me, gave me space to develop, taught me, encouraged me, provided critical feedback, and befriended me.  What more could one ask for?  I never figured out what was ``unnatural" for him.
	
	I feel quite {inpsired} by his investment in my development.   He helped me find a perfect thesis question, and gave me the space and support to solve the problem.  Eventually, with {some good ideas and} blue collar effort, I found a neat solution, at least for binary trees.  A nonuniqueness issue interfered with a certain construction preventing a proof for $d$-ary trees when $d\ge 3$.  But, as I told Tom, it was time to graduate.  I had a candidate answer for $d\ge 3$, but it remained a conjecture.
	
	Tom waited a time after my graduation before inquiring if I was working on proving the conjecture.  When I said no, he asked if he could take up the quest.  I (nervously) agreed.  What did I miss?  Was there a simple mechanism that allowed one to select a ``natural" solution and show that it had the desired nonnegativity property to complete the construction?  Yikes!  My thesis advisor is working on the unsolved part of my thesis.  Yes, it was a little scary.  I wondered if this was standard issue.
	
	Tom, with his great prowess at building from case-by-case analysis and using extensive computation and pattern recognition to generalize, ended up proving the conjecture \cite{MR1813837}.  Fortunately, he needed some clever ideas related to stochastic monotonicity that ultimately played into his interest in negative dependence, to circumvent the nonuniqueness issue.  In typical Tom fashion, his final approach was elegant.  {Ultimately, I am flattered that he found this project interesting enough to pursue.}
	
	Like he did for so many other former students, postdoctoral scholars, and early career researchers, and even though my research program had shifted away from his beloved interacting particle systems, Tom continued to be a mentor, advocate, and trusted advisor throughout my career.  He was not only a great mathematician -- he was a great person.  I am forever grateful to have intersected paths with him in such significant ways on my semi-random walk through life.

\section*{}	
	
	\begin{center}
		{\large\bf Tom Liggett's work on Interacting Particle Systems}
		
		\vspace{2mm}
		
		Rick Durrett
	\end{center}
	
	Tom Liggett received his Ph.D. from Stanford in 1969, writing a thesis with Sam Karlin on ``Weak convergence of conditioned sums of independent random vectors.'' This was not a very hot topic in probability (his thesis has been cited twice) so he was looking for a new topic to work on. At almost the same time Frank Spitzer began the field of interacting particle systems with his 1970 paper ``Interaction of Markov Processes'' \cite{MR268959}. Most of these processes takes place on a set $S$, e.g., $S = \Z^d$, and the state of the system at time $t$ is a function $\xi_t : S \to F$, where $F$ is some finite set. In words, $\xi_t (x)$ gives the state of $x$ at time $t$. In many examples, $F$ is a two point set so the dynamics are described by giving the rates $c(x,\xi)$ at which $x$ changes its state when the configuration is $\xi$. In typical examples on $\Z^d$ the flip rates are determined by the values of $\xi$ in $\{x\} \cup N$, where $N$ is the set of neighbors of $x$.
	
	Spitzer introduced five examples in his paper. By far the most successful was the simple exclusion process in which particles perform independent random walks subject to the restriction that no two particles could occupy the same site. He also considered the simple exclusion with time change which is a model of a lattice gas, the zero-range process in which particles jump at a rate determined by the number of other particles on the same site, and nearest-neighbor models in one dimension. As Benjamin Weiss said in his summary of the paper in Math Reviews: Results are complete in the case of finite $S$. Not much is proved for infinite $S$ but there are many interesting conjectures. 
	
	Two events combined to bring Tom to the field of interacting particle systems. (i)  Chuck Stone showed Tom a copy of Spitzer's 1970 paper saying ``I think you'll find something interesting in this.'' (ii) He had just completed a paper with Mike Crandall on nonlinear semigroups \cite{MR287357}. The rest, as they say, is history. In 1972 Tom wrote a paper proving the existence of interacting particle systems using ideas he had learned working on semigroups \cite{MR309218}. The existence problem is nontrivial when $S$ is infinite because there is no first jump. At about the same time Ted Harris gave a special construction in the finite range case, but Tom's result was elegant and very general. It has covered almost all of the examples that have been investigated in the last fifty years. 
	
	Tom's 1975 paper \cite{MR402985} with Dick Holley introduced the voter model. Here $F$ is a two-point set of possible opinions, e.g., Democrat $(1)$ or Republican $(0)$. Their simple minded voters wake up at times of their personal rate-one Poisson process and adopt the opinion of a neighbor chosen at random. The key to the analysis of this system is a duality between the voter model and a system of coalescing random walks (CRW). Intuitively the particles in the CRW trace back in time, the origins of the opinions of sites at time $t$. The main utility of duality is that it allows us to study the dynamics of the voter model on an infinite set by running the CRW from any finite initial configurations. Using well-known properties of random walk, Holley and Liggett proved that (i) in dimensions $d \le 2$ the system reached consensus, i.e., the probability of seeing two different opinions in a fixed finite box goes to 0 as $t \to\infty$, and (ii) in $d >2$ there is a one-parameter family of stationary distributions $\nu_\rho$ where $\rho$ is the fractions of 1s in equilibrium. The last result highlights some of the difficulties of many interacting particle systems: there is one-parameter family of stationary distributions all of which are mutually singular.
	
	In 1978 Holley and Liggett turned their attention to the contact process, which is perhaps the simplest particle system \cite{MR488379}. Each site is in state 1 = occupied or 0 = vacant. Deaths ($1\to 0$ transitions) occur at a constant rate, while births ($0\to1$) occur at rate $\lambda$ times the number of occupied neighbors. Harris introduced this process in 1974. It is not hard to show that in the one-dimensional nearest neighbor case, if $\lambda< 1$ then the process dies out (reaches the all 0 state) starting from any finite set of 1s. It is much harder to find a number $\Lambda$ so that if $\lambda >\Lambda$ then the system survives (starting from a single 1 there is a positive probability of not dying out). 
	
	The smallest value of $\lambda$ for which survival occurs is the critical value $\lambda_c$. Harris proved that in the one-dimensional nearest neighbor case $\lambda_c<\infty$. With a little work you can prove $\lambda_c < 57$ with Harris' method. In 1978, Holley and Liggett \cite{MR488379} proved $\lambda_c\le  2$. The method is ingenious. If a particle system is attractive, that is the birth rate is an increasing function of the configuration and the death rate is decreasing, then starting from all 1s, the process converges to a limit. If the limit is not a point mass on the all 0s state then the limit is a nontrivial stationary distribution. If it is, then there is no nontrivial stationary distribution. While this dichotomy is useful, it often takes a lot of work to determine which outcome occurs. The brilliant idea of Holley and Liggett was to find an initial condition that increases in distribution, which guarantees that there is a nontrivial stationary distribution. This paper should come with a warning label ``don't try this at home'' because I have never seen anyone else succeed in using this approach. It is not too hard to find a renewal measure that is a good initial condition but then one must show that for every increasing function $f$,  one has that $\mathbf{E} f(\xi_t)$ increases in $t$.
	
	In 1992 Robin Pemantle introduced the contact process on trees $T^d$ in which each vertex has the same degree $d$. This process is interesting because it has two phase transitions. If $\lambda <\lambda_1$ then the process dies out. If $\lambda_1 <\lambda <\lambda_2$ then the process does not die out but the particles wander off to $\infty$ (i.e., any fixed $x\in T^d$ is not occupied infinitely often). If $\lambda >\lambda_2$ then when the process does not die out, any given site is occupied infinitely often. Pemantle used a number of clever arguments to show that if the degree $d \ge 4$, then $\lambda_1 <\lambda_2$.  In 1996 Liggett \cite{MR1421990} settled the final case of $d = 3$ ($d=2$ is just $\Z$) by showing that $\lambda_1 < .605$ and $\lambda_2 > .609$. The proof is not pretty. Basically he uses Pemantle's approach of finding functions of the configuration that are supermartingles, and then automates the process of improving the bounds. Fortunately for the rest of us, at about the same time Alan Stacey, who was a postdoc at UCLA, came up with a soft proof that on any tree with $d \ge 2$ the two thresholds were different.	
	
	There are many variants on the voter model. In the threshold-one voter model, the very insecure voters change their opinion at rate 1 if at least one neighbor has a different opinion. As is often the case in interacting particle systems nothing interesting happens in a one-dimensional nearest neighbor system. In this case, the two opinions cannot coexist. In 1994 Tom proved (under natural conditions on the set of neighbors) that coexistence  is possible in all other cases \cite{MR1288131}. Using some simple comparisons, it is enough to consider the threshold one contact process (in which births occur at vacant sites when at least one neighbor is occupied) and show that when the neighborhood $N = \{-2, -1, 1, 2\}$ and $\lambda=1$, then the system survives. With a little help from his computer he was able to show there was survival when $\lambda= 0.985$.
	
	In addition to performing intricate computations to solve concrete problems Tom has developed powerful theoretical results. In the 1980s probabilists were interested in proving ``shape theorems'' for randomly growing objects. One of the simplest examples is first passage percolation. In this model there are independent and identically distributed positive times assigned to the edges of $\Z^d$. The passage time $t(x,y)$ is the total time on the fastest path from $x$ to $y$, and we are interested in the asymptotic behavior of $W(t) = \{ x : t(0,x) \le t\}$, i.e., the points that can be reached from the origin by time $t$. To study this growing set, one starts with the observation that for fixed $x\in\Z^d$  and positive integers $l, m,$ and $n$ we have $t(lx,mx) + t(mx,nx) \ge t(lx,nx)$, and then one uses Kingman's subadditive ergodic theorem (with some auxiliary arguments) to show that $W(t)/t$ converges to a convex set that has the same symmetries as those of  $\Z^d$  (that leave the origin fixed). 
	
	Shape theorems have been proved for the contact process, biased voter model, branching random walks, and many other models. In these examples the subadditivity condition does not hold when $l > 0$. So, in 1985 Tom developed an improved version of Kingman's result that covers these examples and others such as the right-edge of the one-dimensional contact process \cite{MR806224}. His result has assumptions that are easy to check and a very clever argument for the hardest part of the proof: a lower bound on $\liminf_{n\to\infty} X(0,n)/n$. Here $X(m,n)$ is the subadditive process.
	
	For interacting particle systems like the contact process, it is easy to prove results when $\lambda$ is much smaller than $\lambda_c$ or much larger than $\lambda_c$. In  1984 I invented a ``block construction'' that is, in a sense, a rigorous version of the physics notion of renormalization. Intuitively, as you look at the system on larger and larger space scales, the parameter flows away from the fixed point at $\lambda_c$. The result in mathematical terms is that you can prove results for all $\lambda>\lambda_c$ by proving results for oriented percolation with $p$ close to 1. The bad news is that the states of the sites in the oriented percolation have a finite range of dependence. The dependence is a little annoying, but since $p$ is close to 1, one can prove the necessary result by adapting the usual ``contour argument'' for proving results in this regime. I don't think Tom was ever a big fan of this technique (which I have used repeatedly over the last 35 years), so he, Robert Schonmann, and Alan Stacey, showed \cite{MR1428500} that one could find an independent percolation with almost the same density that was dominated by the $M$-dependent percolation, eliminating the need to continually reprove things for $M$-dependent oriented percolation.
	
	Tom's books have done much to advance the field. His 1977 St. Flour lecture notes introduced the world to the subject. I learned from a preprint of these notes before I came to UCLA. These notes were made obsolete by his 1985 monograph \cite{MR776231}, but for those who are curious, they are available with Spitzer's 1972 notes (in French) and my 1995 notes on the block construction in a Springer book titled {\it Interacting Particle Systems at St. Flour}. For those who are not familiar with this conference series, it has been held every summer in the city in the south of France that appears in its name. Each year it features three speakers giving a ten lecture series. Sadly this year, the 50th conference, which featured Ivan Corwin, Sylvie M\'el\'eard, and Allan Sly was canceled due to the Covid-19 pandemic.
	
	Tom's 1985 book {\it Interacting Particle Systems} (IPS) brought together many of the results that had been proved since 1970. The first three chapters are still worth reading. They show 1. how to construct the process, 2. explain some useful tools such as coupling, monotonicity, and duality; and in 3. specialize to the case of spin systems in which only one site changes its state when a jump occurs. With these basics introduced, the remaining chapters turn to specific examples: 4. the Ising model, 5. the voter model, 6. the contact process, 7. nearest particle systems, 8. the exclusion process, and 9. linear systems in which sites have nonnegative states. The last topic is a creation of Frank Spitzer for his 1979 Wald lectures. In one example there is an inner product type duality between the smoothing process, in which a site is replaced by an average of its neighbors, and the potlatch which is ``an opulent ceremonial feast at which possessions are given away to attendees.'' See Holley and Liggett's 1981 paper \cite{MR608015} for details.
	
	By the time Tom came to write his 1999 book, the field had become too large to cover, so Tom concentrated his attention on the contact process, the voter model, and the simple exclusion process. The book begins with a quick summary of the background from IPS. At the time of the 1985 book, the contact process was well understood in $d=1$ but not $d> 1$. The 1999 book describes the work of Bezuidenhout and Grimmett that solved all of the major open problems in $d>1$. Due to the clarity of the exposition, most people learn this material from Tom's book rather than the original papers. The 1999 book also has a detailed account of the contact process on trees. Tom was the first to prove some of the results, but more importantly synthesized the results in the literature into a comprehensive treatment. 
	
	The voter model chapter contains results for the threshold one and higher thesholds. The chapter on the exclusion process has sections on the relationship to Burger's equation via hydrodynamic limits, the behavior of tagged particles, and results for the system on $\{1, 2,\ldots, N\}$. Each chapter has a set of notes that place the research in context. 
	
	A list of Tom's publications since 2000 can be found on his UCLA web page. There are 48 papers written with a number of co-authors: Richard Arratia, Phillip Bonacich, Julius Borcea,  Maury Bramson, Peter Branden,  Pietro Caputo, Lincoln Chayes, Subroshekhar Ghosh, Geoffrey Grimmett, Alexander Holroyd, Steve Lippman, Tom Mountford,  Thomas Richthammer, Silke Roles, Dan Romik, Richard Rumlet, Rinaldo Schinazi, Jason Schweinsberg, Jeff Steif, Wenping Tang, Balint Toth, Yuan Zhang, and me. Many of these people were his colleagues or postdocs at UCLA, but many others are successful probabilists who live at various places in the US and around the world.
	The 48 papers since 2000 address a diverse set of topics: finitely dependent colorings, cellular automata, random graphs, social network formation, Shapley values for market games, phylogenetic trees, hard-core interactions, and the mysterious: How to find an extra head (in a sequence of heads and tails). I think many of these collaborations arose like mine with Tom and Yuan Zhang did. We had part of the answer and turned to Tom for help with completing the solution.
	
	According to MathSciNet, Tom has 106 papers that have been cited 3341 times by 2333 authors. His 1985 and 1999 books lead the list, followed by his paper \cite{MR287357} on nonlinear semigroup with Mike Crandall. But even the papers with citations in single digits are interesting and well written.

\section*{}

	\begin{center}
		{\large\bf Tom Liggett's work on Finite Dependence and Colorings}
		
		\vspace{2mm}
		
		Alexander E. Holroyd
	\end{center}
	
	Every mathematician has their style.  One of the thrills of the subject is learning to appreciate this, and experiencing how those styles can enhance and complement each other.  I want to explain an astonishing and wonderful facet of Tom Liggett's style that was perhaps literally unique -- a superpower, one might say.  He once explained it to me, with the candor, directness and self-awareness that many of us so cherished, something like this:
	
	\begin{quotation}\em``When I approach a problem, I typically have no idea at all how to solve it.  But I do something that, so far as I know, no-one else does: I play around by hand with small cases, $n=1$, $n=2$, and so on, and just look for patterns.''
	\end{quotation}
	
	It seems to me that Tom's self-deprecating analysis wasn't quite right. The method he described is I think familiar to almost all mathematicians.  The difference is simply that Tom was thousands of times better at it than anyone else!
	
	I will try to illustrate this by sharing the story of our collaboration on finite dependence, which turned out to be one of my favourite pieces of mathematics.
	A stochastic process $X=(X_v)_{v\in V}$ indexed by a metric space $V$ is called {\it $k$-dependent} if, for any two sets $A,B\subseteq V$ at distance greater than $k$ from each other, the two random vectors $(X_v)_{v\in A}$ and $(X_v)_{v\in B}$ are independent.  A process is {\it finitely dependent} if it is $k$-dependent for some finite $k$.  Finite dependence is arguably the strongest and simplest mixing condition, and has applications in statistical physics, computer science and probability theory.  My own introduction to the concept (before I knew Tom and the other authors personally) was the extremely useful paper \cite{MR1428500} on stochastic domination.  (It has hundreds of citations.  I believe I used the main result seven separate times in my PhD thesis.)
	
	Despite the simplicity of the definition, finite dependence is not so easy to understand, even in the setting of stationary processes on the integer line $V=\Z$ (the focus henceforth).  The obvious way to construct a stationary finitely dependent process $X$ is as a {\it block factor}: all the randomness comes from an i.i.d.\ family $U=(U_v)_{v\in \Z}$, and $X$ is given by $X_v=f(U_{v+1},U_{v+2},\ldots, U_{v+r})$ where $f$ is a fixed function of $r$ arguments.  Are these the only examples?  This question seems to have its origins in the 1960s.  It was open for several decades until the first counterexamples were found in the 1990s -- see e.g.\ \cite{MR1245304}.  A series of subsequent papers explored the idea further, but all known counterexamples had the feel of artificial processes constructed specifically for the purpose.  A consensus emerged that perhaps the only `natural' stationary finitely dependent processes are block factors.
	
	The main question I want to discuss arose in conversations between myself, Itai Benjamini, Oded Schramm and Benjamin Weiss in 2008.  Motivated in part by distributed computing and statistical physics, we were interested in the interaction between mixing properties and hard local constraints on a stochastic process.  The canonical constraint is (proper) coloring: a process $X=(X_v)_{v\in\Z}$ is a {\it $q$-coloring} if each $X_v$ takes values in the finite set of `colors' $\{1,\ldots,q\}$, and almost surely $X_v\neq X_{v+1}$ for all $v$.  We quickly realized that there was a very simple question we couldn't answer:
	\begin{quotation}\centering\em Is there a stationary finitely dependent coloring of $\Z$?\end{quotation}
	
	We all became convinced that the answer must be no, due in part to several negative results.  We established that if such a process existed, it could not be a block factor (leaving only the `unnatural' finitely dependent non-block-factors); it could not be Markov or hidden Markov (with a finite underlying state space); it could not be a $1$-dependent $3$-coloring (the first non-trivial case).
	
	Oded Schramm became particularly interested in the problem (he told me he was `obsessed' with it).  Among other things, he translated it into a question about existence of a certain exotic Hilbert space, which was passed around to various experts; he introduced a Fourier transform approach aimed at ruling out all finitely dependent $3$-colorings.  After Oded's tragic death in September 2008 I continued to work on the problem, sometimes in consultation with others, but I gradually became convinced that it was extremely difficult, and started to lose hope of ever knowing the answer.
	
	That feeling changed immediately when I mentioned the problem to Tom in 2011, and, to my delight, he told me that it was exactly the kind of question he liked.  (This was by no means a given. Tom was incredibly disciplined and selective in his focus.  If it was not his kind of problem, he would decide quickly and say so.)   Progress was initially rapid and exciting.  After a week or so, Tom had found a new elegant proof of the impossibility of $1$-dependent $3$-coloring.  A week or two later he told me he thought he had proved the impossibility of $1$-dependent $4$-coloring, and that he would start on $1$-dependent $5$-coloring next.  A few weeks after that he told me there was a mistake in the last proof, and he now believed that there \emph{was} a $1$-dependent $4$-coloring!  (Given my prior experience with the problem, I was skeptical).  In fact, a few days later, Tom produced an extraordinary explicit formula for the cylinder probabilities of a putative stationary $1$-dependent $4$-coloring.
	
	Here is the formula (after some superficial simplifications from me).
	\begin{equation}\label{formula}
	P\binom{y}{z}=2^{-m}\sum_{w\in \text{DD}(m-1)}(-1)^{|w|}\, c(w,y,z)\, \mu(y_w).
	\end{equation}
	Here, the $4$ colors are identified with the binary column vectors of signs, $1,2,3,4=\binom++,\binom+-,\binom-+,\binom--$.   The left side is the cylinder probability of a proper coloring $x=(x_1,\ldots,x_n)=\binom y z$ of length $n$, written as a $2$-by-$n$ matrix of signs with rows $y,z\in\{+,-\}^n$.  The variable $m$ is the number of runs of $+$s and $-$s in $y$.  The sum is over the set $\text{DD}(m-1)$ of \emph{dispersed Dyck words} of length $m-1$, i.e.\ elements of $\{\circ,\langle\,,\rangle\}^{m-1}$ consisting of a concatenation of $\circ$s and Dyck words (Dyck words are legal bracket-sequences of equal numbers of $\langle$s and $\rangle$s); we think of the symbols of $w$ as aligned with the internal run-boundaries of $y$, in order; $|w|$ is the number of $\langle$s in $w$.  The coefficient $c(w,y,z)$ is the product over the internal run-boundaries of $y$ of the $z$-symbol immediately left or right of the boundary, according to whether the corresponding $w$-symbol is $\langle$ or $\rangle$ respectively, or neither if it is $\circ$.  The string $y_w$ is $y$ with some of its internal runs sign-flipped, in such a way that run-boundaries corresponding to $\langle$s and $\rangle$s are eliminated.  For example, if $y=\p\p\m\m\p\p\p\m$ and $w=\circ\langle\,\rangle$ then $y_w=\p\p\m\m\m\m\m\m$, because the run $\p\p\p$ is flipped to $\m\m\m$.  Finally, $\mu(y_w)$ is the probability that a uniformly random permutation of length $n+1$ has its descents in exactly the positions of the $-$s of $y_w$.
	
	How was this obtained?  Playing around with small cases, looking for patterns!  Tom told me he could probably never fully explain the procedure to anyone else's satisfaction, but the idea as I understand it is as follows.  The requirement of a stationary $1$-dependent $4$-coloring imposes various constraints (equalities and inequalities) on the cylinder probabilities.  For instance, $P(132)+P(142)=P(1{*}2)=P(1)P(2)$.  On their face, the constraints seem insufficient to uniquely determine the probabilities, but when there is a choice, one can try to pick the simplest or most parsimonious possibility, and then work through the consequences.  Tom apparently played this delicate game, gradually refining the choices, and eventually guessing an incredibly complicated pattern.  I am doubtful whether anyone else could ever have achieved this.
	
	We were able to check that the formula satisfied all the requirements, with one important exception.  We were unable to prove that $P(\cdot)$ was non-negative, as required for a probability.  There followed three years of attempts by Tom and me to prove this.  We checked billions of cases by computer.  We constructed more and more elaborate arguments involving induction, inclusion-exclusion and  transforms that proved many sub-cases; we formulated further conjectures about partial sums that also seemed to be non-negative; there were tantalizing connections with completely unexpected entries from the Online Encyclopedia of Integer Sequences; hundreds of pages of \LaTeX\ were exchanged.  In another twist, the formula for $P(x)$ appeared to be symmetric under the action of all permutations of the $4$ colors.  Even symmetry between $y$ and $z$ is not at all clear, although Tom had expected it from his construction. We were not able to prove this either.  We knew what the marginal distributions for any choice of one or two colors had to be: the appearances of one color are equal in distribution to the locations of $HT$ in a sequence of fair coin flips, while any two colors combined correspond to the descents in a sequence of i.i.d.\ uniform variables.  I asked every probabilist who would listen whether they could think of \emph{any} process with these properties.  I got plenty of fascinating suggestions, but no solution.  After three years, we were getting weary, and wondering whether to simply give up and publish \eqref{formula} as a conjecture.
	
	The eventual breakthrough came in 2014, when we found that the $P$ of \eqref{formula} appeared to satisfy the remarkably simple but unusual recursion
	\begin{equation}\label{recursion}
	P(x)=\begin{cases}\frac{1}{2n+2} \sum_{i=1}^n P(\widehat{x}_i) &\text{if $x$ is a proper coloring};\\
	0&\text{otherwise.}
	\end{cases}
	\end{equation}
	where $x$ has length $n$, and $\widehat{x}_i=x_1\cdots x_{i-1}x_{i+1} \cdots x_n$ denotes $x$ with the $i$th element deleted.  This was found initially by comparing \eqref{formula} with recursions satisfied by $\mu$, and then checking cases by computer.  Once one has \eqref{recursion}, it is not too difficult to prove that it is indeed satisfied by the $P$ of \eqref{formula}, and of course non-negativity then follows immediately.  In fact one can simply take \eqref{recursion} as the definition of $P$, and prove directly that the resulting process is $1$-dependent.  This proof is short but mysterious - we still do not have good intuition for why it works.  The solution to \eqref{recursion} is clearly symmetric under permutations of the colors, so the same applies to \eqref{formula}.
	
	So, the answer to the 2008 question is yes!  The critical number of colors is $4$ -- there is a stationary $1$-dependent $4$-coloring but no $3$-coloring.  Moreover, coloring distinguishes between finitely dependent processes and block factors, providing a very satisfying answer to this much older question (as well as some other questions from the finite dependence literature).  In fact, one can deduce that any non-trivial local constraint system (suitably defined) distinguishes between the two classes of processes.  And Oded's exotic Hilbert space exists.
	
	Formally, the final proof does not need Tom's miraculous formula \eqref{formula}, because one can instead work entirely from the recursion \eqref{recursion}.  But there is no way I could have found \eqref{recursion} without \eqref{formula} -- despite its simplicity, the recursion seems so arbitrary, and so unrelated to the problem it solves.
	
	There is one more miracle.  For which $(k,q)$ does a stationary $k$-dependent $q$-coloring exist?  Not $(1,3)$, but $(1,4)$.  We can do $(3,3)$ as well: starting from our $1$-dependent $4$-coloring, we can eliminate every $4$ by replacing it with the smallest-numbered color not among its two neighbors, to get a $3$-dependent $3$-coloring. The question is trivially monotone in $k$ and $q$, and the cases $k=0$ and $q=2$ are easily ruled out, so that only leaves $(2,3)$.  Tom and I were almost ready to submit our paper, stating that existence of a $2$-dependent $3$-coloring seemed to be a very difficult open question.  We decided to check, for completeness, that just applying the same recursion \eqref{recursion} would not work.  Well, it turns out that it does work.  More precisely, one can apply \eqref{recursion} to construct a stationary $q$-coloring for any $q\geq 2$ (with a different normalizing constant in place of $1/(2n+2)$ to maintain a probability measure).  With $q=4$ the resulting process is $1$-dependent.  With $q=3$ it is $2$-dependent.   For all other values of $q$, the process is not finitely dependent.  Only $q=3$ and $q=4$ work, and they give precisely the minimal possible cases.  (One startling consequence is that conditioning the $1$-dependent $4$-coloring to have no $4$s gives a $2$-dependent $3$-coloring.) There is surely something fundamental and mysterious at work here.
	
	Everything discussed above appears in \cite{MR3570073}.  Subsequent papers (e.g.\ \cite{MR4079439}) have explored and extended the ideas, although many mysteries remain.  There is still essentially only one known construction of a finitely dependent coloring (or even a finitely dependent process satisfying non-trivial hard constraints) -- all others are variations on the same theme.  It is unknown whether automorphism-invariant finitely dependent colorings exist on regular trees or general Euclidean lattices.  Motivated in part by the method that led to the discovery of the formula \eqref{formula}, we conjecture that the $1$-dependent $4$-coloring is unique, but currently there is little idea how to prove this.  Perhaps most perplexing and fascinating of all, no-one really understands why the construction works, even though the proof is quite short, and later discoveries have helped to broaden the context somewhat.
	
	Tom was one of the finest people I have known.  He was an incredible mathematician, mentor and friend, and the best collaborator one could hope for.  I will miss him.  I will miss his no-nonsense advice whenever I face a difficult decision. I will miss his brilliance and his mathematical superpower whenever I encounter a beautiful question.

\section*{}

	\begin{center}
		{\large\bf Tom Liggett and Aldous' Spectral Gap Conjecture}
		
		\vspace{2mm}
		
		Pietro Caputo
	\end{center}
	
	\subsection*{The problem}
	Consider a connected, undirected graph $G=(V,E)$ with $n$ vertices. A particle configuration is a permutation $\sigma$ assigning a label to each vertex, with the interpretation that the particle with label $\sigma(i)$ occupies vertex $i$. 
	The interchange process is the continuous-time Markov chain on the symmetric group $S_n$ of all $n!$ permutations defined as follows: for each edge $(i,j)\in E$ independently, at the arrival times of a rate 1 Poisson process, interchange the labels at vertex $i$ and vertex $j$.
	
	The interchange process is reversible, irreducible, and the stationary distribution, call it  $\nu$, is uniform on $S_n$. A classical measure of the speed of convergence to stationarity is the spectral gap $\lambda_{IP}(G)$, that is the smallest non-zero eigenvalue of the transition rate matrix. If we follow the evolution of a single particle, the projected dynamics coincides with the continuous-time simple random walk on $G$. This is a reversible irreducible process,  
	and hence has a positive spectral gap $\lambda_{RW}(G)$. The contraction principle shows that $\lambda_{IP}(G) \leq \lambda_{RW}(G)$.
	Aldous' conjecture stated that 
	\begin{equation}\label{eq:conj}
	\lambda_{IP}(G) = \lambda_{RW}(G),
	\end{equation}
	for all graphs $G$.
	The conjecture appeared on David Aldous' webpage and it was mentioned as an open problem in a well-known online monograph by Aldous and Fill, {\it Reversible Markov Chains and Random Walks on Graphs}. 
	
	The exclusion process with $k\leq n$ particles on $G$ is obtained by painting black the particles with labels $1,\dots,k$ and white the particles with labels $k+1,\dots,n$, so that configurations are now given by the $\binom{n}{k}$ possible positions of $k$ indistinguishable particles on $n$ vertices with at most one particle per vertex.   A consequence of the the identity \eqref{eq:conj} is that the spectral gap of the exclusion process with $k$ particles on $G$ is the same for all  $k=1,\dots,n-1$. At the time of its formulation in 1992, the conjecture was known to hold only for complete graphs and for star graphs, where the special symmetries allow the computation of the spectrum of the interchange process via Fourier analysis on the symmetric group (see the work of Diaconis and Shahshahani \cite{MR626813} and a subsequent work by Flatto, Odlyzko, and Wales). A natural generalisation suggests that the conjectured identity \eqref{eq:conj} should extend to any weighted graph, that is when  the Poisson clock at the edge $(i,j)$ has rate $c(i,j)$, where $c(i,j)\geq 0$ are arbitrary edge weights or conductances.
	
	\subsection*{Background}
	Understanding the behaviour of the spectral gap of the exclusion process and other interacting particle systems has been a fundamental problem for many years, with motivations ranging from statistical physics to theoretical computer science. 
	Starting with pioneering work of D.\ Aldous and P.\ Diaconis on random walks on groups in the 80's the analysis of convergence to stationarity became a central theme in the probability literature. 
	Spectral gap estimates played a key role in the ground-breaking works of Jerrum, Sincalir, Holley, Stroock, Varadhan, H.T. Yau, Quastel, Martinelli and many others on Glauber and Kawasaki dynamics for lattice spins systems. These works developed efficient recursive techniques that allowed in certain cases the determination of the correct scaling of the spectral gap with the system size.     
	For instance with these methods one could easily prove that the spectral gap of the interchange process on cubic boxes of $\bbZ^d$ scales diffusively with the linear size, that is $\lambda_{IP}(G)\asymp \ell^{-2}$ if $G$ is a lattice cube with side $\ell$, which is indeed the scaling of the spectral gap $\lambda_{RW}(G)$ of the lattice Laplacian for any $d\geq 1$. In particular, for certain special families of graphs, with these methods one could establish the asymptotic equivalence $ \lambda_{IP}(G) \asymp \lambda_{RW}(G)$. 
	The more  fundamental and algebraic flavour of the conjecture \eqref{eq:conj} however seemed to call for new methods and new ideas. Some progress was obtained for special families of graphs. The conjecture was shown to hold for all weighted trees by Handjani and Jungreis \cite{MR1419872}, who discovered a neat recursive argument.  A version of the conjecture in terms of energy levels of one-dimensional quantum spin Hamiltonians was proven by Koma and Nachtergaele. For the case of large lattice boxes the identity \eqref{eq:conj} was shown to hold asymptotically as the side of the box diverges by Conomos, Starr and Morris. On the other hand a fully algebraic approach was developed by Cesi  who showed the validity of \eqref{eq:conj} for all complete multipartite graphs. 
	
	\subsection*{Working with Tom on the conjecture} 
	As most graduate students in this field in the 90's I learned about interacting particle systems from Tom Liggett's book {\it Interacting Particle Systems} \cite{MR776231}. At that time his book was so influential among us that we often referred to it as {\em the bible}. 
	But I had never met Tom before he invited me for a one-year visit at UCLA in 2008. That's also the first time  I met with Thomas Richthammer, then a postdoc in probability at UCLA.
	After the first semester was over we decided it was time to work on some project together, the three of us. 
	When Thomas and I proposed to work on Aldous' conjecture, Tom looked at us  and said ``this is a very hard problem, you know we're not going to solve it''. He 
	loved the problem, he had been playing with it, and had already made a couple of serious attempts at it in the past. He told us he had convinced himself he would never succeed in solving it.  But the smile on his face was clearly saying he was still intrigued by the challenge. He seemed happy to have a chance to think about that again.
	
	I had some ideas in the back of my head that I wanted to try, so I started playing with various kinds of recursions, but none of that seemed to  go anywhere. 
	Tom was more to the point. He said: let's write a proof for arbitrary weighted graphs with $n=4$ vertices. He did that, in all details. It was quite elementary in the end and very specific to $n=4$, but it inspired us,  we discussed it over and over. At the same time I was looking for a simple interpretation of the recursive technique that allowed the authors of \cite{MR1419872} to solve the problem for arbitrary trees, one that could be modified in order to attack more general graphs. I like to reformulate their argument as follows. Let $G$ be a tree with $n$ vertices and let $\cE_G(f)$ denote the Dirichlet form of the interchange process on $G$, so that $\lambda_{IP}(G)$ is the largest real number $\lambda$ such that  
	\begin{equation}\label{eq:lam1}
	\cE_G(f)\geq \lambda \,\var_n(f)\,,
	\end{equation}
	for any function $f:S_n\mapsto\bbR$, where $\var_n(f)$ denotes the variance of $f$ under the uniform distribution $\nu$. Using orthogonality of eigenfunctions  
	one can reduce the proof of \eqref{eq:conj} to showing that \eqref{eq:lam1} holds with $\lambda=\lambda_{RW}(G)$ for all $f$ satisfying 
	\begin{equation}\label{eq:var}
	\var_n(f) = \nu[\var_n(f | \sigma(i))],
	\end{equation}
	where $\var_n(f | \sigma(i))$ denotes the variance of $f$ with respect to the conditional distribution $\nu(\cdot | \sigma(i))$ obtained by revealing the label at vertex $i$, and $\nu[\,\cdot\,]$ denotes the expectation with respect to $\nu$. Here $i$ is an arbitrary fixed vertex. If we choose $i$ to be a leaf of the tree $G$, then removing it together with the edge incident to it one obtains a new tree $G_i$ with $n-1$ vertices for which, inductively, we may assume the identity \eqref{eq:conj} to be true, that is 
	\begin{equation}\label{eq:lam2}
	\cE_{G_i}(g)\geq \lambda_{RW}(G_i) \,\var_{n-1}(g)\,,
	\end{equation}
	for all $g:S_{n-1}\mapsto\bbR$. Since $\var_n(f | \sigma(i)) = \var_{n-1}(g)$ where $g$ is the function obtained from $f$ by freezing the label at vertex $i$ as $\sigma(i)$,  \eqref{eq:var} and \eqref{eq:lam2} imply
	\begin{align}\label{eq:lam3}
	&\lambda_{RW}(G_i)\var_n(f) =  \lambda_{RW}(G_i)\,\nu\left[\var_{n-1}(g)\right] \nonumber \\
	&\leq \nu\left[\cE_{G_i}(g)\right] \leq \cE_G(f)\,,
	\end{align}
	where, in the last step, we have used an obvious monotonicity associated to the removal of $i$ and the edge incident to it, namely that
	\begin{equation}\label{eq:lam4}
	\nu\left[\cE_{G_i}(g)\right] \leq \cE_G(f)\,.
	\end{equation}
	Thus, the identity \eqref{eq:conj} is proved if one shows that 
	\begin{equation}\label{eq:lam5}
	\lambda_{RW}(G_i)\geq \lambda_{RW}(G).
	\end{equation}
	The inequality \eqref{eq:lam5} is not as obvious as \eqref{eq:lam4} but it is not difficult to prove. A way to interpret \eqref{eq:lam5}
	is by viewing the random walk on $G_i$ as the trace on $V\setminus\{i\}$ of  the random walk on $G$, so that \eqref{eq:lam5} follows from a suitable contraction. This ended the proof for trees.
	
	Naive attempts at extending this strategy of proof to graphs with cycles ran into problems. For instance, if $G_i$ is obtained by removing a vertex in a cycle together with all edges incident to it, then \eqref{eq:lam4} is again obviously true but there is no hope that \eqref{eq:lam5} continues to hold. 
	
	The three of us discussed during lunch breaks at the UCLA faculty lounge, and then tried to combine this point of view with the details of various specific examples we could  by then solve explicitly.  It took several weeks until  the first major breakthrough finally arrived: we should devise a {\em different} transformation $G\mapsto G_i$ which maps a graph with vertex set $V$ to a graph with vertex set $V\setminus\{i\}$ in such a way that both \eqref{eq:lam4} and \eqref{eq:lam5} were satisfied, and the transformation should be obtained by redistributing the edge weights in a way that would be natural from the point of view of electric network theory, that is  when the edge weights are interpreted as conductances. Namely, if $G$ has vertex set $V$ and edge weights $c(\cdot,\cdot)$, the graph $G_i$ should be defined as the weighted graph with vertex set  $V\setminus\{i\}$ and new weights $c'(\cdot,\cdot)$ given by:
	\begin{equation}\label{eq:rates}
	c'(j,k) = c(j,k) + \frac{c(i,j)c(i,k)}
	{\sum_{\ell} c(i,\ell)}\,, \qquad j,k \in V\setminus\{i\}.
	\end{equation}
	For instance, for 3-vertex  networks this is the familiar series reduction of electrical circuits. In $4$-vertex networks  this is known as the $Y-\Delta$ or  star-triangle transformation.  
	A very pleasant feature of this transformation is that it is precisely the one for which the random walk on $G_i$ can be seen as the trace on $V\setminus\{i\}$ of the random walk on $G$. In particular, one can check that \eqref{eq:lam5} is always satisfied. If one could prove that \eqref{eq:lam4} is also satisfied then the proof would be complete. We made several tests and everything indicated that this was the right recursive strategy. We got very excited and thought that if such an inequality must hold for arbitrary choices of the weights then a proof should be simple. It was not.

	For almost two months we mostly worked on our own, and met two or three times a week for short briefings. During our discussions we drew many pictures of edges emanating from a vertex but the inequality \eqref{eq:lam4} kept escaping us. The edges became tentacles and we started calling it the {\em octopus inequality}. A name that Tom enjoyed very much. The inequality can be reformulated in more explicit terms as follows: 
	\begin{align}\label{eq:octopus}
	&\sum_{\ell} c(i,\ell)\, \nu\left[(f\circ\tau_{i\ell} - f)^2\right] \nonumber\\
	&\geq \frac12\sum_{j,k}\frac{c(i,k)c(i,j)}{\sum_{\ell} c(i,\ell)} \, \nu\left[(f\circ\tau_{jk} - f)^2\right],
	\end{align}
	where $\tau_{jk}$ denotes the transposition at the edge $(j,k)$. The inequality should hold for every fixed $i\in V$, for all functions $f$ on $S_n$ and for all choices of non-negative edge weights $c(\cdot,\cdot)$. In abstract terms, the inequality says that a certain $n!\times n!$ matrix $C$ is positive semidefinite. The matrix $C$ formally looks like a transition rate matrix except that it involves both positive and negative rates. We agreed at some point that the octopus had the flavour of a correlation inequality but could not immediately use this analogy to prove its validity.
	
	One day Tom told us more details about his own attempts at catching the octopus. For small values of $n$ he had found a very useful way of rewriting $C$ as a covariance matrix of some simple random variables, which implied the desired positivity. This was the last  breakthrough we needed. Thomas quickly understood what had to be tried to reproduce a similar argument for higher values of $n$ and 
	a few more days of hard work finally allowed us to finish the proof. 
	We were running out of time since both Thomas and I were getting ready to return  to Europe. Moreover, in the meantime we had learned that our recursive approach based on electric network reduction had just been independently discovered by Dieker \cite{MR2600660}, so we were no longer alone in our hunt for the octopus.  
	The proof of the octopus inequality turned out to be a somewhat mysterious combination of subtle manipulations combined with some elementary group-theoretic facts such as properties of  cosets
	of subgroups of even and odd permutations in $S_n$. 
	
	\subsection*{Perspective}
	The result was received with enthusiasm and got published in the prestigious Journal of the AMS \cite{MR2629990}. Later, Tom several times referred to this as one of his best works. 
	The overall strategy of proof is indeed very pleasant and it brought to light a new functional inequality, the octopus inequality, which is of interest in its own right. 
	The proof of the octopus inequality however did not illuminate much on its meaning, and it seemed natural to hope for an alternative, perhaps simpler proof.  An algebraic approach to the octopus, with some simplifications and some new insights has been given recently by Cesi. However, after more than ten years we still do not have a substantially  more appealing proof of it. 
	
	The octopus inequality has been recently used as a tool to obtain new comparison inequalities in works of Chen, Alon, Kozma, Hermon, and Salez. It turned out that it is sufficiently powerful to determine the correct scaling of the mixing time of the interchange process for certain families of graphs. 
	
	The last time I met with Tom we talked about a conjecture I had since our joint work, that the identity \eqref{eq:conj} should extend to a larger class of continuous-time Markov chains on $S_n$, obtained by replacing graphs with hypergraphs. More precisely, assign to each $A\subset V$  a weight $\alpha_A\geq 0$ and for each $A$ independently, perform perfect shuffles of the particles in $A$ at the arrival times of a Poisson process with rate $\alpha_A$.  Call this the $\alpha$-shuffle process. Following just one particle reveals a random walk on $V$ with transition rates $$c_\alpha(i,j) = \sum_{A\subset V:\,A\ni i,j}\frac{\alpha_A}{|A|},$$
	where $|A|$ is the cardinality of the subset $A$.  I conjectured
	that the $\alpha$-shuffle process and the random walk with rates $c_\alpha(\cdot,\cdot)$ have the same spectral gap. When the weights $\alpha$ are such that $\alpha_A=0$ unless $|A|=2$ this coincides with \eqref{eq:conj}. Indeed, in this case the interchange process with rates $c(i,j)=\alpha_{\{i,j\}}$ coincides with the $\alpha$-shuffle process run at twice the speed.  One hope is to discover a larger octopus-like inequality that would shed some new light on the original octopus. However, a direct approach based on an analogue of the one-vertex reduction seems to fail.  
	Tom liked this conjecture and strongly encouraged me to write about it.
	
	\subsection*{Postscript by David Aldous}
	
	(a) People sometimes ask ``how did the spectral gap conjecture arise?'', and it was merely that we couldn't think
	of any function of the process that might relax more slowly than the single particle process.
	
	(b) It is often asserted that in analysis, there are deep inequalities, but no deep equalities.
	I tell students that ``general identities'' are almost always
	best understood as the same quantity calculated in two different ways (e.g. the Parseval relation)
	or the same diagram interpreted in two different ways (e.g. duality in Interacting Particle Systems via graphical representations).
	But even with the simplified approach to the octopus, the argument still seems deep.
	So is this result an exception to the general assertion, or have we all missed a simpler proof?
	
	(c) It is ultimately a result about symmetric matrices, and of course symmetric matrices arise
	in many fields within the mathematical sciences.  I had vaguely hoped that the result would  
	relate to problems in some other field, but so far no such connection has appeared.

\bibliographystyle{alpha}

\begin{thebibliography}{BGM93}
	
	\bibitem[BBL09]{MR2476782}
	Julius Borcea, Petter Br\"{a}nd\'{e}n, and Thomas~M. Liggett.
	\newblock Negative dependence and the geometry of polynomials.
	\newblock {\em J. Amer. Math. Soc.}, 22(2):521--567, 2009.
	
	\bibitem[BGM93]{MR1245304}
	Robert~M. Burton, Marc Goulet, and Ronald Meester.
	\newblock On {$1$}-dependent processes and {$k$}-block factors.
	\newblock {\em Ann. Probab.}, 21(4):2157--2168, 1993.
	
	\bibitem[CL71]{MR287357}
	M.~G. Crandall and T.~M. Liggett.
	\newblock Generation of semi-groups of nonlinear transformations on general
	{B}anach spaces.
	\newblock {\em Amer. J. Math.}, 93:265--298, 1971.
	
	\bibitem[CLR10]{MR2629990}
	Pietro Caputo, Thomas~M. Liggett, and Thomas Richthammer.
	\newblock Proof of {A}ldous' spectral gap conjecture.
	\newblock {\em J. Amer. Math. Soc.}, 23(3):831--851, 2010.
	
	\bibitem[Die10]{MR2600660}
	A.~B. Dieker.
	\newblock Interlacings for random walks on weighted graphs and the interchange
	process.
	\newblock {\em SIAM J. Discrete Math.}, 24(1):191--206, 2010.
	
	\bibitem[DS81]{MR626813}
	Persi Diaconis and Mehrdad Shahshahani.
	\newblock Generating a random permutation with random transpositions.
	\newblock {\em Z. Wahrsch. Verw. Gebiete}, 57(2):159--179, 1981.
	
	\bibitem[HHL20]{MR4079439}
	Alexander~E. Holroyd, Tom Hutchcroft, and Avi Levy.
	\newblock Mallows permutations and finite dependence.
	\newblock {\em Ann. Probab.}, 48(1):343--379, 2020.
	
	\bibitem[HJ96]{MR1419872}
	Shirin Handjani and Douglas Jungreis.
	\newblock Rate of convergence for shuffling cards by transpositions.
	\newblock {\em J. Theoret. Probab.}, 9(4):983--993, 1996.
	
	\bibitem[HL75]{MR402985}
	Richard~A. Holley and Thomas~M. Liggett.
	\newblock Ergodic theorems for weakly interacting infinite systems and the
	voter model.
	\newblock {\em Ann. Probability}, 3(4):643--663, 1975.
	
	\bibitem[HL78]{MR488379}
	R.~Holley and T.~M. Liggett.
	\newblock The survival of contact processes.
	\newblock {\em Ann. Probability}, 6(2):198--206, 1978.
	
	\bibitem[HL81]{MR608015}
	Richard Holley and Thomas~M. Liggett.
	\newblock Generalized potlatch and smoothing processes.
	\newblock {\em Z. Wahrsch. Verw. Gebiete}, 55(2):165--195, 1981.
	
	\bibitem[HL16]{MR3570073}
	Alexander~E. Holroyd and Thomas~M. Liggett.
	\newblock Finitely dependent coloring.
	\newblock {\em Forum Math. Pi}, 4:e9, 43, 2016.
	
	\bibitem[Lig72]{MR309218}
	Thomas~M. Liggett.
	\newblock Existence theorems for infinite particle systems.
	\newblock {\em Trans. Amer. Math. Soc.}, 165:471--481, 1972.
	
	\bibitem[Lig85a]{MR806224}
	Thomas~M. Liggett.
	\newblock An improved subadditive ergodic theorem.
	\newblock {\em Ann. Probab.}, 13(4):1279--1285, 1985.
	
	\bibitem[Lig85b]{MR776231}
	Thomas~M. Liggett.
	\newblock {\em Interacting particle systems}, volume 276 of {\em Grundlehren
		der Mathematischen Wissenschaften [Fundamental Principles of Mathematical
		Sciences]}.
	\newblock Springer-Verlag, New York, 1985.
	
	\bibitem[Lig94]{MR1288131}
	Thomas~M. Liggett.
	\newblock Coexistence in threshold voter models.
	\newblock {\em Ann. Probab.}, 22(2):764--802, 1994.
	
	\bibitem[Lig96]{MR1421990}
	Thomas~M. Liggett.
	\newblock Branching random walks and contact processes on homogeneous trees.
	\newblock {\em Probab. Theory Related Fields}, 106(4):495--519, 1996.
	
	\bibitem[Lig00]{MR1813837}
	Thomas~M. Liggett.
	\newblock Monotonicity of conditional distributions and growth models on trees.
	\newblock {\em Ann. Probab.}, 28(4):1645--1665, 2000.
	
	\bibitem[LSS97]{MR1428500}
	T.~M. Liggett, R.~H. Schonmann, and A.~M. Stacey.
	\newblock Domination by product measures.
	\newblock {\em Ann. Probab.}, 25(1):71--95, 1997.
	
	\bibitem[Spi70]{MR268959}
	Frank Spitzer.
	\newblock Interaction of {M}arkov processes.
	\newblock {\em Advances in Math.}, 5:246--290 (1970), 1970.
	
\end{thebibliography}

\end{document}